\documentclass[reqno]{amsart}
\usepackage{amsfonts}
\usepackage{graphics,graphicx,amssymb}
\usepackage{amsmath,amsthm}
\usepackage{hyperref}
\numberwithin{equation}{section}
\newtheorem{theorem}{Theorem}[section]
\newtheorem{corollary}{Corollary}[section]
\newtheorem{proposition}{Proposition}[section]
\newtheorem{lemma}{Lemma}[section]
\newtheorem{definition}{Definition}[section]

\newfont{\got}{eufm9 scaled 1095}

\newfont{\w}{msbm9 scaled\magstep1}

\begin{document}

\title{ALMOST COMPLEX CONNECTIONS ON ALMOST COMPLEX MANIFOLDS WITH NORDEN METRIC}
\author{MARTA TEOFILOVA}
\date{}

\maketitle

\begin{abstract} A four-parametric family of linear connections preserving the almost
complex structure is defined on an almost complex manifold with
Norden metric. Necessary and sufficient conditions for these
connections to be natural are obtained. A two-parametric family of
complex connections is studied on a conformal K\"{a}hler manifold
with Norden metric. The curvature tensors of these connections are
proved to coincide.

\noindent \emph{2010 Mathematics Subject Classification}: 53C15,
53C50, 32Q60.

\noindent \emph{Keywords}: Almost complex manifold with Norden
metric, almost complex connection.

\end{abstract}

\maketitle

\section*{Introduction}
Almost complex manifolds with Norden metric were first studied by
A.~P.~Norden \cite{N}. These manifolds are introduced in
\cite{Gri-Mek-Dje} as generalized $B$-manifolds. A classification of
the considered manifolds with respect to the covariant derivative of
the almost complex structure is obtained in \cite{Gan} and two
equivalent classifications are given in \cite{Gan-Gri-Mih2,Gan-Mih}.

An important problem in the geometry of almost complex manifolds
with Norden metric is the study of linear connections preserving the
almost complex structure or preserving both, the structure and the
Norden metric. The first ones are called \emph{almost complex}
connections, and the second ones are called \emph{natural}
connections. A special type of a natural connection is the canonical
one. In \cite{Gan-Mih} it is proved that on an almost complex
manifold with Norden metric there exists a unique canonical
connection. The canonical connection and its conformal group on a
conformal K\"{a}hler manifold with Norden metric are studied in
\cite{Gan-Gri-Mih2}.

In \cite{Teo1} we have studied the Yano connection on a complex
manifold with Norden metric and in \cite{Teo2} we have proved that
the curvature tensors of the canonical connection and the Yano
connection coincide on a conformal K\"{a}hler manifold with Norden
metric.

In the present paper we define a four-parametric family of almost
complex connections on an almost complex manifold with Norden
metric. We find necessary and sufficient conditions for these
connections to be natural. By this way we obtain a two-parametric
family of natural connections on an almost complex manifold with
Norden metric. We study a two-parametric family of complex
connections on a conformal K\"{a}hler manifold with Norden metric,
obtain the form of the K\"{a}hler curvature tensor corresponding to
each of these connections and prove that these tensors coincide.

\section{Preliminaries}

Let $(M,J,g)$ be a $2n$-dimensional almost complex manifold with
Norden metric, i.e. $J$ is an almost complex structure and $g$ is a
metric on $M$ such that
\begin{equation}
J^{2}X=-X,\qquad g(JX,JY)=-g(X,Y)  \label{1.1}
\end{equation}
for all differentiable vector fields $X$, $Y$ on $M$, i.e. $X,Y\in %
\mbox{\got X}(M)$.

The associated metric $\widetilde{g}$ of $g$, given by $\widetilde{g}%
(X,Y)=g(X,JY)$, is a Norden metric, too. Both metrics are
necessarily neutral, i.e. of signature $(n,n)$.

Further, $X,Y,Z,W$ ($x,y,z,w$, respectively) will stand for
arbitrary differentiable vector fields on $M$ (vectors in $T_{p}M$,
$p\in M$, respectively).

If $\nabla $ is the Levi-Civita connection of the metric $g$, the
tensor field $F$ of type $(0,3)$ on $M$ is defined by
$F(X,Y,Z)=g\left( (\nabla _{X}J)Y,Z\right)$ and has the following
symmetries
\begin{equation}
F(X,Y,Z)=F(X,Z,Y)=F(X,JY,JZ).  \label{1.3}
\end{equation}

Let $\left\{ e_{i}\right\} $ ($i=1,2,\ldots ,2n$) be an arbitrary
basis of $ T_{p}M$ at a point $p$ of $M$. The components of the
inverse matrix of $g$ are denoted by $g^{ij}$ with respect to the
basis $\left\{ e_{i}\right\} $. The Lie 1-forms $\theta $ and
$\theta^{\ast}$ associated with $F$, and the Lie vector $\Omega$,
corresponding to $\theta$, are defined by, respectively
\begin{equation}
\theta (z)=g^{ij}F(e_{i},e_{j},z), \qquad \theta^{\ast}=\theta \circ
J, \qquad \theta (z)=g(z,\Omega ). \label{theta}
\end{equation}

The Nijenhuis tensor field $N$ for $J$ is given by
$N(X,Y)=[JX,JY]-[X,Y]-J[JX,Y]-J[X,JY]$. The corresponding tensor of
type (0,3) is given by $N(X,Y,Z)=g(N(X,Y),Z)$. In terms of $\nabla
J$ this tensor is expressed in the following way
\begin{equation}
N(X,Y)=(\nabla_{X}J)JY - (\nabla_{Y}J)JX + (\nabla_{JX}J)Y -
(\nabla_{JY}J)X. \label{N}
\end{equation}
It is known \cite{N-N} that the almost complex structure is complex
if and only if it is integrable, i.e. $N=0$. The associated tensor
$\widetilde{N}$ of $N$ is defined by \cite{Gan}
\begin{equation}
\widetilde{N}(X,Y)=(\nabla_{X}J)JY + (\nabla_{Y}J)JX +
(\nabla_{JX}J)Y + (\nabla_{JY}J)X, \label{Ntilde}
\end{equation}
and the corresponding tensor of type (0,3) is given by
$\widetilde{N}(X,Y,Z)=g(\widetilde{N}(X,Y),Z)$.

A classification of the almost complex manifolds with Norden metric
is introduced in \cite{Gan}, where eight classes of these manifolds
are characterized according to the properties of $F$. The three
basic classes $\mathcal{W}_{i}$ ($i=1,2,3$) and the class
$\mathcal{W}_{1}\oplus\mathcal{W}_{2}$ are given by

$\bullet$ the class $\mathcal{W}_{1}$:
\begin{equation}\label{w1}
\begin{array}{l}
F(X,Y,Z)=\frac{1}{2n}\left[ g(X,Y)\theta (Z)+g(X,JY)\theta
(JZ)\right. \smallskip \\
\qquad \qquad \qquad \qquad \quad \left. +g(X,Z)\theta
(Y)+g(X,JZ)\theta (JY)\right];
\end{array}
\end{equation}

$\bullet$ the class $\mathcal{W}_{2}$ of \emph{the special complex
manifolds with Norden metric}:
\begin{equation}\label{w2}
F(X,Y,JZ)+F(Y,Z,JX)+F(Z,X,JY)=0,\quad \theta =0;
\end{equation}

$\bullet$ the class $\mathcal{W}_{3}$ of \emph{the quasi-K\"{a}hler
manifolds with Norden metric}:
\begin{equation}\label{w3}
F(X,Y,Z)+F(Y,Z,X)+F(Z,X,Y)=0 \Leftrightarrow \widetilde{N}=0;
\end{equation}

$\bullet$ the class $\mathcal{W}_{1}\oplus\mathcal{W}_{2}$ of
\emph{the complex manifolds with Norden metric}:
\begin{equation}\label{w1w2}
F(X,Y,JZ)+F(Y,Z,JX)+F(Z,X,JY)=0 \Leftrightarrow  N=0.
\end{equation}
The special class $\mathcal{W}_{0}$ of \emph{the K\"{a}hler
manifolds with Norden metric} is characterized by $F=0$.

A $\mathcal{W}_{1}$-manifold with closed Lie 1-forms $\theta$ and
$\theta^{\ast}$ is called \emph{a conformal K\"{a}hler manifold with
Norden metric}.

Let $R$ be the curvature tensor of $\nabla $, i.e.
$R(X,Y)Z=\nabla _{X}\nabla _{Y}Z-\nabla _{Y}\nabla _{X}Z-\nabla _{\left[ X,Y%
\right] }Z$ and $R(X,Y,Z,W)=g\left( R(X,Y)Z,W\right)$.

A tensor $L$ of type (0,4) is said to be \emph{curvature-like} if it
has the properties of $R$, i.e.
$L(X,Y,Z,W)=-L(Y,X,Z,W)=-L(X,Y,W,Z)$,
$L(X,Y,Z,W)+L(Y,Z,X,W)+L(Z,X,Y,W)=0$. Then, the Ricci tensor
$\rho(L)$ and the scalar curvatures $\tau(L)$ and $ \tau^{\ast}(L)$
of $L$ are defined by:
\begin{equation}
\begin{array}{c}
\rho(L)(y,z)=g^{ij}L(e_{i},y,z,e_{j}),\smallskip\\\tau(L)
=g^{ij}\rho(L)(e_{i},e_{j}),\quad \tau^{\ast}(L)=g^{ij}\rho(L)
(e_{i},Je_{j}).
\end{array}
\label{Ricci, tao}
\end{equation}
A curvature-like tensor $L$ is called \emph{a K\"{a}hler tensor} if
$L(X,Y,JZ,JW) = - L(X,Y,Z,W)$.

Let $S$ be a tensor of type (0,2). We consider the following tensors
\cite{Gan-Gri-Mih2}:
\begin{equation}\label{psi}
\begin{array}{l}
\psi_{1}(S)(X,Y,Z,W)=g(Y,Z)S(X,W)-g(X,Z)S(Y,W)\smallskip\\
\qquad \qquad \qquad \qquad\hspace{0.17in} + \hspace{0.045in}
g(X,W)S(Y,Z) - g(Y,W)S(X,Z),\medskip\\
\psi_{2}(S)(X,Y,Z,W)=g(Y,JZ)S(X,JW) - g(X,JZ)S(Y,JW)\smallskip\\
\qquad \qquad \qquad \qquad \hspace{0.17in}+ \hspace{0.045in}
g(X,JW)S(Y,JZ) - g(Y,JW)S(X,JZ), \medskip\\
\pi_{1}=\frac{1}{2}\psi_{1}(g),\qquad\quad
\pi_{2}=\frac{1}{2}\psi_{2}(g),\qquad\quad
\pi_{3}=-\psi_{1}(\widetilde{g})=\psi_{2}(\widetilde{g}).
\end{array}
\end{equation}
The tensor $\psi_{1}(S)$ is curvature-like if $S$ is symmetric, and
the tensor $\psi_{2}(S)$ is curvature-like is $S$ is symmetric and
hybrid with respect to $J$, i.e. $S(X,JY)=S(Y,JX)$. The tensors
$\pi_{1} - \pi_{2}$ and $\pi_{3}$ are K\"{a}hlerian.

\section{Almost complex connections on almost complex manifolds
with Norden metric}

In this section we study almost complex connections and natural
connections on almost complex manifolds with Norden metric. First,
let us recall the following
\begin{definition}\cite{Ko-No}\label{def-complex} A linear connection $\nabla^{\prime}$ on
an almost complex manifold $(M,J)$ is said to be \emph{almost
complex} if $\nabla^{\prime}J=0$.
\end{definition}
\begin{theorem}\label{th2.1} On an almost complex manifold with Norden metric
there exists a 4-parametric family of almost complex connections
$\nabla^{\prime}$ with torsion tensor $T$ defined by, respectively:
\begin{equation}\label{complex connections}
\begin{array}{l}
g\big( \nabla^{\prime}_{X}Y - \nabla_{X}Y, Z \big) =
\frac{1}{2}F(X,JY,Z)+ t_{1}\big\{ F(Y,X,Z)\smallskip\\+ F(JY,JX,Z)
\big\} +t_{2}\big\{ F(Y,JX,Z) - F(JY,X,Z) \big\}\smallskip\\ +
t_{3}\big\{ F(Z,X,Y)+F(JZ,JX,Y) \big\}\smallskip\\ + t_{4}\big\{
F(Z,JX,Y) - F(JZ,X,Y) \big\},
\end{array}
\end{equation}
\begin{equation}\label{torsion}
\begin{array}{l}
\quad T(X,Y,Z) =t_{1}\big\{F(Y,X,Z) - F(X,Y,Z)+
F(JY,JX,Z)\smallskip\\ \quad - F(JX,JY,Z)\big\}+
\big(\frac{1}{2}-t_{2}\big)\big\{F(X,JY,Z) - F(Y,JX,Z)\big\}
\smallskip\\ \quad+ t_{2}\big\{F(JX,Y,Z) - F(JY,X,Z)\big\}
+ 2t_{3}F(JZ,JX,Y)\smallskip\\ \quad + 2t_{4}F(Z,JX,Y),
\end{array}
\end{equation}
where $t_{i}\in \mathbb{R},\hspace{0.1in} i=1,2,3,4$.
\begin{proof}
By (\ref{complex connections}), (\ref{1.3}) and direct computation,
we prove that $\nabla^{\prime}J=0$, i.e. the connections
$\nabla^{\prime}$ are almost complex.
\end{proof}
\end{theorem}

By (\ref{w3}) and (\ref{w1w2}) we obtain the form of the almost
complex connections $\nabla^{\prime}$ on the manifolds belonging to
the classes $\mathcal{W}_{1}\oplus \mathcal{W}_{2}$ and
$\mathcal{W}_{3}$ as follows, respectively
\begin{corollary}\label{cor-1} On a complex manifold
with Norden metric there exists a 2-parametric family of complex
connections $\nabla^{\prime}$ defined by
\begin{equation}\label{complexconn-1}
\begin{array}{l}
\nabla^{\prime}_{X}Y = \nabla_{X}Y + \frac{1}{2}(\nabla_{X}J)JY +
p \big\{(\nabla_{Y}J)X + (\nabla_{JY}J)JX \big\} \smallskip\\
\qquad \quad + q \big\{(\nabla_{Y}J)JX - (\nabla_{JY}J)X \big\},
\end{array}
\end{equation}
where $p=t_{1}+t_{3}$, $q=t_{2}+t_{4}$.
\end{corollary}
\begin{corollary}\label{cor-2} On a quasi-K\"{a}hler manifold with
Norden metric there exists a 2-parametric family of almost complex
connections $\nabla^{\prime}$ defined by
\begin{equation}\label{complexconn-2}
\begin{array}{l}
\nabla^{\prime}_{X}Y = \nabla_{X}Y + \frac{1}{2}(\nabla_{X}J)JY +
s \big\{(\nabla_{Y}J)X + (\nabla_{JY}J)JX \big\} \medskip\\
\qquad \quad + t \big\{(\nabla_{Y}J)JX - (\nabla_{JY}J)X \big\},
\end{array}
\end{equation}
where $s=t_{1}-t_{3}$, $t=t_{2}-t_{4}$.
\end{corollary}
\begin{definition}\cite{Gan-Mih}\label{def-natural} A linear connection $\nabla^{\prime}$ on an
almost complex manifold with Norden metric $(M,J,g)$ is said to be
\emph{natural} if
\begin{equation}
\nabla^{\prime}J = \nabla^{\prime}g =0 \quad(\Leftrightarrow
\nabla^{\prime}g=\nabla^{\prime}\widetilde{g}=0).
\end{equation}
\end{definition}
\begin{lemma}\label{lema}
Let $(M,J,g)$ be an almost complex manifold with Norden metric and
let $\nabla^{\prime}$ be an arbitrary almost complex connection
defined by (\ref{complex connections}). Then
\begin{equation}\label{nabla g}
\begin{array}{l}
\big(\nabla^{\prime}_{X}g\big)(Y,Z) = (t_{2} +
t_{4})\widetilde{N}(Y,Z,X) - (t_{1} +
t_{3})\widetilde{N}(Y,Z,JX),\smallskip\\
\big(\nabla^{\prime}_{X}\widetilde{g}\big)(Y,Z) = -(t_{1} +
t_{3})\widetilde{N}(Y,Z,X) - (t_{2} + t_{4})\widetilde{N}(Y,Z,JX).
\end{array}
\end{equation}
\end{lemma}
Then, by help of Theorem \ref{th2.1} and Lemma \ref{lema} we prove
\begin{theorem}\label{th2.2}
An almost complex connection $\nabla^{\prime}$ defined by
(\ref{complex connections}) is natural on an almost complex manifold
with Norden metric if and only if $t_{1}=-t_{3}$ and $t_{2}=-t_{4}$,
i.e.
\begin{equation}\label{natural}
\begin{array}{l}
g\big(\nabla^{\prime}_{X}Y - \nabla_{X}Y, Z \big) =
\frac{1}{2}F(X,JY,Z) \smallskip\\ \qquad\qquad \qquad \qquad
\hspace{0.1in} + t_{1}N(Y,Z,JX) - t_{2}N(Y,Z,X).
\end{array}
\end{equation}
\end{theorem}
The equation (\ref{natural}) defines a 2-parametric family of
natural connections on an almost complex manifold with Norden metric
and non-integrable almost complex structure. In particular, by
(\ref{w3}) and (\ref{nabla g}) for the manifolds in the class
$\mathcal{W}_{3}$ we obtain
\begin{corollary}\label{cor natural}
Let $(M,J,g)$ be a quasi-K\"{a}hler manifold with Norden metric.
Then, the connection $\nabla^{\prime}$ defined by
(\ref{complexconn-2}) is natural for all $s,t \in \mathbb{R}$.
\end{corollary}

If $(M,J,g)$ is a complex manifold with Norden metric, then from
(\ref{w1w2}) and (\ref{natural}) it follows that there exists a
unique natural connection $\nabla^{\prime}$ in the family
(\ref{complexconn-1}) which has the form
\begin{equation}\label{D}
\begin{array}{l}
\nabla^{\prime}_{X}Y = \nabla_{X}Y + \frac{1}{2}(\nabla_{X}J)JY.
\end{array}
\end{equation}

\begin{definition}\cite{Gan-Mih}\label{def-can} A natural
connection $\nabla^{\prime}$ with torsion tensor $T$ on an almost
complex manifold with Norden metric is said to be \emph{canonical}
if
\begin{equation}\label{canonical-cond}
T(X,Y,Z) + T(Y,Z,X) - T(JX,Y,JZ) - T(Y,JZ,JX) = 0.
\end{equation}
\end{definition}

Then, by applying the last condition to the torsion tensors of the
natural connections (\ref{natural}), we obtain
\begin{proposition}Let $(M,J,g)$ be an almost complex manifold with
Norden metric. A natural connection $\nabla^{\prime}$ defined by
(\ref{natural}) is canonical if and only if $t_{1}=0$,
$t_{2}=\frac{1}{8}$. In this case (\ref{natural}) takes the form
\begin{equation*}\label{can-connection}
\begin{array}{l}
2g\big(\nabla^{\prime}_{X}Y - \nabla_{X}Y, Z\big) = F(X,JY,Z) -
\frac{1}{4}N(Y,Z,X).
\end{array}
\end{equation*}
\end{proposition}

Let us remark that G. Ganchev and V. Mihova \cite{Gan-Mih} have
proven that on an almost complex manifold with Norden metric there
exists a unique canonical connection. The canonical connection of a
complex manifold with Norden metric has the form (\ref{D}).

Next, we study the properties of the torsion tensors (\ref{torsion})
of the almost complex connections $\nabla^{\prime}$.

The torsion tensor $T$ of an arbitrary linear connection is said to
be \emph{totally antisymmetric} if $T(X,Y,Z)=g(T(X,Y),Z)$ is a
3-form. The last condition is equivalent to
\begin{equation}\label{3-form}
T(X,Y,Z)=-T(X,Z,Y).
\end{equation}
Then, having in mind (\ref{torsion}) we obtain that the torsion
tensors of the almost complex connections $\nabla^{\prime}$ defined
by (\ref{complex connections}) satisfy the condition (\ref{3-form})
if and only if $t_{1}=t_{2}=t_{3}=0$, $t_{4}=\frac{1}{4}$. Hence, we
prove the following

\begin{theorem}
Let $(M,J,g)$ be an almost complex manifold with Norden metric and
non-integrable almost complex structure. Then, on $M$ there exists a
unique almost complex connection $\nabla^{\prime}$ in the family
(\ref{complex connections}) whose torsion tensor is a 3-form. This
connection is defined by
\begin{equation} \label{anti-sym}
\begin{array}{l}
g\big(\nabla^{\prime}_{X}Y -
\nabla_{X}Y,Z\big)=\frac{1}{4}\big\{2F(X,JY,Z)+
F(Z,JX,Y)\smallskip\\ \qquad \qquad \qquad \qquad \hspace{0.1in}-
F(JZ,X,Y) \big\}.
\end{array}
\end{equation}
\end{theorem}

By Corollary \ref{cor-2} and the last theorem we obtain
\begin{corollary}
On a quasi-K\"{a}hler manifold with Norden metric there exists a
unique natural connection $\nabla^{\prime}$ in the family
(\ref{complexconn-2}) whose torsion tensor is a 3-form. This
connection is given by
\begin{equation}\label{bismut}
\begin{array}{l}
\nabla^{\prime}_{X}Y = \nabla_{X}Y +
\frac{1}{4}\big\{2(\nabla_{X}J)JY - (\nabla_{Y}J)JX +
(\nabla_{JY}J)X \big\}.
\end{array}
\end{equation}
\end{corollary}

Let us remark that the connection (\ref{bismut}) can be considered
as an analogue of the Bismut connection \cite{Bi}, \cite{Gau} in the
geometry of the almost complex manifolds with Norden metric.

Let us consider symmetric almost complex connections in the family
(\ref{complex connections}). By (\ref{torsion}) and (\ref{N}) we
obtain
\begin{equation}\label{T-condition}
\begin{array}{l}
T(X,Y) - T(JX,JY) = \frac{1}{2}N(X,Y).
\end{array}
\end{equation}
From (\ref{T-condition}), (\ref{torsion}) and (\ref{N}) it follows
that $T=0$ if and only if $N=0$ and $t_{1}=t_{3}=t_{4}=0$,
$t_{2}=\frac{1}{4}$. Then, it is valid the following
\begin{theorem}
Let $(M,J,g)$ be a complex manifold with Norden metric. Then, on $M$
there exists a unique complex symmetric connection $\nabla^{\prime}$
belonging to the family (\ref{complexconn-1}) which is given by
\begin{equation}\label{Yano}
\begin{array}{l}
\nabla^{\prime}_{X}Y=\nabla_{X}Y+\frac{1}{4}\big\{(\nabla_{X}J)JY +
2(\nabla_{Y}J)JX - (\nabla_{JX}J)Y\big\}.
\end{array}
\end{equation}
\end{theorem}
The connection (\ref{Yano}) is known as the Yano connection
\cite{Ya1,Ya2}.

We give a summery of the obtained results for the 4-parametric
family of almost complex connections $\nabla^{\prime}$ in the
following
\begin{center}
$\underset{\text{Table 1}}{
\begin{tabular}{l|c|c|c|}
\cline{2-4} & \multicolumn{3}{|c|}{\small{Class manifolds}} \\
\hline \multicolumn{1}{|l|}{\small{Connection type}} &
\small{$\mathcal{W}_{1}\oplus\mathcal{W}_{2}\oplus\mathcal{W}_{3}$}
& \small{$\mathcal{W}_{1}\oplus \mathcal{W}_{2}$}
& \small{$\mathcal{W}_{3}$} \\
\hline \multicolumn{1}{|l|}{\small{almost complex}} &
\small{$t_{1},t_{2},t_{3},t_{4}\in\mathbb{R}$} & \small{$p,q \in
\mathbb{R}$}
& \small{$s,t \in \mathbb{R}$} \\
\hline \multicolumn{1}{|l|}{\small{natural}} &
\small{$t_{1}=-t_{3}$, $t_{2}=-t_{4}$} & \small{$p=q=0$}
& \small{$s,t \in \mathbb{R}$}  \\
\hline \multicolumn{1}{|l|}{\small{canonical}} &
\small{$t_{1}=t_{3}=0$, $t_{2}=-t_{4}=\frac{1}{8}$} &
\small{$p=q=0$}
& \small{$s=0$, $t=\frac{1}{4}$} \\
\hline \multicolumn{1}{|l|}{\small{$T$ is a 3-form}} &
\small{$t_{1}=t_{2}=t_{3}=0$, $t_{4}=\frac{1}{4}$} &
\small{$\nexists$}
& \small{$s=0$, $t=-\frac{1}{4}$} \\
\hline \multicolumn{1}{|l|}{\small{symmetric}} & \small{$\nexists$}
& \small{$p=0$, $q=\frac{1}{4}$} & \small{$\nexists$}
\\ \hline
\end{tabular}}$
\end{center}

\section{Complex connections on conformal K\"{a}hler manifolds with Norden metric}

Let $(M,J,g)$ be a $\mathcal{W}_{1}$-manifold with Norden metric and
consider the 2-parametric family of complex connections
$\nabla^{\prime}$ defined by (\ref{complexconn-1}). By (\ref{w1}) we
obtain the form of $\nabla^{\prime}$ on a $\mathcal{W}_{1}$-manifold
as follows
\begin{equation}\label{Q}
\begin{array}{l}
\nabla^{\prime}_{X}Y=\nabla_{X}Y + \frac{1}{4n}\big\{g(X,JY)\Omega -
g(X,Y)J\Omega + \theta(JY)X \smallskip\\  -\theta(Y)JX \big\}+
\frac{p}{n}\big\{\theta(X)Y+\theta(JX)JY\big\} +
\frac{q}{n}\big\{\theta(JX)Y - \theta(X)JY\big\}.
\end{array}
\end{equation}
Then, by (\ref{Q}) and straightforward computation we prove
\begin{theorem}
Let $(M,J,g)$ be a conformal K\"{a}hler manifold with Norden metric
and $\nabla^{\prime}$ be an arbitrary complex connection in the
family (\ref{complexconn-1}). Then, the K\"{a}hler curvature tensor
$R^{\prime}$ of $\nabla^{\prime}$ has the form
\begin{equation*} \label{Rw10}
R^{\prime} = R -\frac{1}{4n}\big\{\psi_{1} + \psi_{2}\big\}(S) -
\frac{1}{8n^{2}}\psi_{1}(P) -
\frac{\theta(\Omega)}{16n^{2}}\big\{3\pi_{1}+\pi_{2}\big\} +
\frac{\theta(J\Omega)}{16n^{2}}\pi_{3},
\end{equation*}
where $S$ and $P$ are defined by, respectively:
\begin{equation}\label{SM}
\begin{array}{l}
S(X,Y) = \big(\nabla_{X}\theta\big)JY +
\frac{1}{4n}\big\{\theta(X)\theta(Y) -
\theta(JX)\theta(JY)\big\},\smallskip\\
P(X,Y) = \theta(X)\theta(Y) + \theta(JX)\theta(JY).
\end{array}
\end{equation}
\end{theorem}

By (\ref{Q}) we prove the following
\begin{lemma}\label{lemma w1}
Let $(M,J,g)$ be a $\mathcal{W}_{1}$-manifold and $\nabla^{\prime}$
be an arbitrary complex connection in the family
(\ref{complexconn-1}). Then, the covariant derivatives of $g$ and
$\widetilde{g}$ are given by
\begin{equation}\label{gw1}
\begin{array}{l}
\big(\nabla^{\prime}_{X}g\big)(Y,Z) =-\frac{2}{n}\big\{[p
\theta(X)+q \theta(JX)]g(Y,Z) \smallskip\\ \qquad\qquad\qquad\qquad
+ [p
\theta(JX)-q \theta(X)]g(Y,JZ)\big\},\smallskip\\
\big(\nabla^{\prime}_{X}\widetilde{g}\big)(Y,Z) =
\frac{2}{n}\big\{[p \theta(JX)-q \theta(X)]g(Y,Z) \smallskip\\
\qquad\qquad\qquad\qquad - [p \theta(X)+q \theta(JX)]g(Y,JZ)\big\}.
\end{array}
\end{equation}
\end{lemma}

It is well-known \cite{Ko-No} that the curvature tensor $R^{\prime}$
and the torsion tensor $T$ of an arbitrary linear connection
$\nabla^{\prime}$ satisfy the second Bianchi identity, i.e.
\begin{equation}\label{BianchiII}
\underset{X,Y,Z}{\mathfrak{S}}\big\{
\big(\nabla^{\prime}_{X}R^{\prime}\big)(Y,Z,W) +
R^{\prime}\big(T(X,Y),Z,W\big)\big\} =0,
\end{equation}
where $\mathfrak{S}$ is the cyclic sum over $X,Y,Z$.

From (\ref{Q}) it follows that the torsion tensor of an arbitrary
connection $\nabla^{\prime}$ in the family (\ref{complexconn-1}) has
the following form on a $\mathcal{W}_{1}$-manifold
\begin{equation}\label{torsionw1}
\begin{array}{l}
T(X,Y)=\frac{1-4q}{4n}\big\{\theta(X)JY -\theta(Y)JX -\theta(JX)Y
+ \theta(JY)X \big\}\smallskip\\
\qquad\qquad\quad\hspace{0.06in}
+\hspace{0.03in}\frac{p}{n}\big\{\theta(X)Y - \theta(Y)X +
\theta(JX)JY - \theta(JY)JX \big\}
\end{array}
\end{equation}

Let us denote $\tau^{\prime}=\tau(R^{\prime})$ and $\tau^{\prime
\ast}=\tau^{\ast}(R^{\prime})$. We establish the following

\begin{theorem}
Let $(M,J,g)$ be a conformal K\"{a}hler manifold with Norden metric,
and $\tau^{\prime}$ and $\tau^{\prime \ast}$ be the scalar
curvatures of the K\"{a}hler curvature tensor $R^{\prime}$
corresponding to the complex connection $\nabla^{\prime}$ defined by
(\ref{complexconn-1}). Then, the function $\tau^{\prime} +
i\tau^{\prime \ast}$ is holomorphic on $M$ and the Lie 1-forms
$\theta$ and $\theta^{\ast}$ are defined in a unique way by
$\tau^{\prime}$ and $\tau^{\prime \ast}$ as follows\emph{:}
\begin{equation}\label{Lie forms}
\theta = 2n\hspace{0.014in} d\big(\arctan
\frac{\tau^{\prime}}{\tau^{\prime \ast}}\big), \qquad \theta^{\ast}
= -2n\hspace{0.014in} d \big(\ln\sqrt{\tau^{\prime \hspace{0.01in}
2} + \tau^{\prime \ast \hspace{0.01in} 2}}\hspace{0.02in} \big).
\end{equation}
\begin{proof}
By (\ref{BianchiII}) and (\ref{torsionw1}) we obtain
\begin{equation}\label{Rzaw1}
\begin{array}{c}
\big(\nabla^{\prime}_{X}R^{\prime}\big)(Y,Z,W) +
\big(\nabla^{\prime}_{Y}R^{\prime}\big)(Z,X,W) +
\big(\nabla^{\prime}_{Z}R^{\prime}\big)(X,Y,W)\medskip\\
=\frac{4q-1}{2n}\big\{\theta(X)R^{\prime}(JY,Z,W)-
\theta(JX)R^{\prime}(Y,Z,W)\medskip\\-\theta(Y)R^{\prime}(JX,Z,W) +
\theta(JY)R^{\prime}(X,Z,W) \medskip\\ + \theta(Z)R^{\prime}(JX,Y,W)
-\theta(JZ)R^{\prime}(X,Y,W) \big\}
\medskip\\ -\frac{2p}{n}\big\{\theta(X)R^{\prime}(Y,Z,W)  +
\theta(JX)R^{\prime}(JY,Z,W)  - \theta(Y)R^{\prime}(X,Z,W)
\medskip\\- \theta(JY)R^{\prime}(JX,Z,W) +
\theta(Z)R^{\prime}(X,Y,W)
 + \theta(JZ)R^{\prime}(JX,Y,W) \big\}.
\end{array}
\end{equation}
Then, having in mind the equalities (\ref{gw1}) and their analogous
equalities for $g^{ij}$, we find the total traces of the both sides
of (\ref{Rzaw1}) and get
\begin{equation}\label{tau-ast}
\begin{array}{l}
d\tau^{\prime} = \frac{1}{2n}\big\{\tau^{\prime \ast}\theta -
\tau^{\prime}\theta^{\ast}\big\},\qquad  d\tau^{\prime \ast} = -
\frac{1}{2n}\big\{\tau^{\prime}\theta + \tau^{\prime
\ast}\theta^{\ast}\big\}.
\end{array}
\end{equation}
From (\ref{tau-ast}) it follows immediately that $d\tau^{\prime
\ast} \circ J = -d\tau^{\prime}$, i.e. the function $\tau^{\prime} +
i\tau^{\prime \ast}$ is holomorphic on $M$ and the equalities
(\ref{Lie forms}) hold.
\end{proof}
\end{theorem}

\bigskip

\noindent\emph{Marta Teofilova\\
Department of Geometry\\
Faculty of Mathematics and Informatics\\
University of Plovdiv\\
236 Bulgaria Blvd.\\
4003 Plovdiv, Bulgaria\\
e-mail:} \verb"marta@uni-plovdiv.bg"

\end{document}